\documentclass[10pt,a4paper]{article}

\usepackage{amsthm}
\usepackage{amssymb}
\usepackage{latexsym}
\usepackage{amsmath}
\usepackage{euscript}
\usepackage{graphicx}
\usepackage{fancyhdr}
\usepackage{calc}
\usepackage{romannum}

\newtheorem{Def}{Definition}[section]
\newtheorem{Lem}[Def]{Lemma}

\newtheorem{Thm}[Def]{Theorem}
\newtheorem{Prop}[Def]{Proposition}

\def\trace{{\mathop{\rm trace}}}
\def\realpart{{\mathop{\rm Re}}}
\def\imagpart{{\mathop{\rm Im}}}
\def\Vol{{\mathop{\rm Vol}}}

\begin{document}
\pagenumbering{arabic}

\title{Generalized Lagrangian mean curvature flow in K\"ahler manifolds that are almost Einstein}
\author{Tapio Behrndt}
\maketitle

\begin{abstract}
We introduce the notion of K\"ahler manifolds that are almost Einstein and we define a generalized mean curvature vector field along submanifolds in them. We prove that Lagrangian submanifolds remain Lagrangian, when deformed in direction of the generalized mean curvature vector field. For a K\"ahler manifold that is almost Einstein, and which in addition has a trivial canonical bundle, we show that the generalized mean curvature vector field of a Lagrangian submanifold is the dual vector field associated to the Lagrangian angle.
\end{abstract}

\section{Introduction}
In a Calabi-Yau manifold with parallel holomorphic volume form $\Omega$ there is a distinguished class of submanifolds called special Lagrangian submanifolds. These are oriented Lagrangian submanifolds, which are calibrated with respect to $\realpart\;\Omega$. Special Lagrangian submanifolds have received a lot of attention since the work by Strominger, Yau and Zaslow \cite{SYZ}, where mirror symmetry is related to special Lagrangian torus fibrations.

The notion of special Lagrangian submanifolds can be generalized to the case when the ambient manifold is almost Calabi-Yau. An almost Calabi-Yau manifold is a K\"ahler manifold together with a non-vanishing, not necessarily parallel, holomorphic volume form. A nice property of almost Calabi-Yau manifolds is that they appear in infinite dimensional families, i.e. the moduli space of almost Calabi-Yau structures is infinite dimensional, while Calabi-Yau structures only appear in finite dimensional families due to the theorem of Tian and Todorov \cite{Tian}, \cite{Todorov} and Yau's proof of the Calabi conjecture \cite{Yau}. Choosing a generic almost Calabi-Yau metric is therefore a much more powerful thing to do than choosing a generic Calabi-Yau metric and, as in the study of moduli spaces of $J$-holomorphic curves, this could be of importance for the study of moduli spaces of special Lagrangian submanifolds as conjectured by Joyce in \cite{Joyce2}. Another nice feature of almost Calabi-Yau manifolds is that explicit almost Calabi-Yau metrics on compact manifolds are known, while there are no non-trivial Calabi-Yau metrics on compact manifolds explicitly known. For instance a quintic in $\mathbb{CP}^4$ equipped with the restriction of the Fubini-Study metric is an almost Calabi-Yau manifold.

Special Lagrangian submanifolds in (almost) Calabi-Yau manifolds have been studied extensively by many authors but up to date there is no general method known how to construct examples of special Lagrangian submanifolds. However, since special Lagrangian submanifolds are calibrated submanifolds they are volume minimizing in their homology class and one is tempted to construct special Lagrangian submanifolds by mean curvature flow of Lagrangian submanifolds. The existence of the Lagrangian mean curvature flow in K\"ahler-Einstein manifolds was first proved by Smoczyk \cite{Smoczyk}. Smoczyk shows that the mean curvature flow of a given compact Lagrangian submanifold remains Lagrangian as long as the flow exists. Thus the problem is to find conditions such that the Lagrangian mean curvature flow exists for all time and converges to a special Lagrangian submanifold. One attempt to this was done by Thomas and Yau \cite{ThomasYau}, where they conjecture that a Lagrangian submanifold satisfying a certain stability condition converges smoothly by Lagrangian mean curvature flow to a non-singular special Lagrangian submanifold in the same homology class. In general there are two problems occurring. Firstly one expects that the evolving Lagrangian submanifold develops a finite time singularity. There are only a few longtime convergence results known for Lagrangian mean curvature flow, for instance by Smoczyk \cite{Smoczyk1}, Smoczyk and Wang \cite{SmoczykWang} and Wang \cite{Wang}. The second problem which occurs is that there exist Lagrangian submanifolds without regular Lagrangian volume minimizers in their homology classes. Examples of such Lagrangian submanifolds were found by Wolfson in \cite{Wolfson}.

In this paper we introduce the notion of K\"ahler manifolds that are almost Einstein (in particular, these contain the class of almost Calabi-Yau manifolds), and we define a generalized mean curvature vector field along submanifolds in them. We show that Lagrangian submanifolds remain Lagrangian under deformation in direction of the generalized mean curvature vector field and we obtain a generalized version of Smoczyk's result. Therefore we call the deformation of Lagrangian submanifolds in direction of the generalized mean curvature vector field a generalized Lagrangian mean curvature flow. We show that the generalized Lagrangian mean curvature flow is the negative gradient flow of the volume functional of some conformally rescaled metric. Moreover, if the ambient manifold is almost Calabi-Yau, then we prove that the one-form associated to the generalized mean curvature vector field of a Lagrangian submanifold is the differential of the Lagrangian angle. As a consequence we show that if the initial Lagrangian has zero Maslov class, then the generalized Lagrangian mean curvature flow can be integrated to a scalar equation.

We remark here that recently, after the first version of the present paper, Smoczyk and Wang showed that in every almost K\"ahler manifold that admits an Einstein connection there exists a generalized mean curvature vector field with the property that Lagrangian submanifolds remain Lagrangian under the deformation in its direction. The generalized Lagrangian mean curvature flow introduced by them contains ours in K\"ahler manifolds that are almost Einstein as an example (see \cite{SmoczykWang2} for more details).

The author would like to thank his supervisor Dominic Joyce for pointing out this problem to him, for discussions about it, and for several corrections. This work was supported by a Sloane Robinson Graduate Award of the Lincoln College, by a scholarship of the British Chamber of Commerce in Germany, and by an EPSRC Research Studentship.

\section{Lagrangian mean curvature flow in K\"ahler-Einstein manifolds}

We first recall the definition of the mean curvature flow. Let $M$ be a Riemannian manifold and let $N$ be a submanifold of $M$ given by an immersion $F_0:N\rightarrow M$. Throughout this paper the term submanifold will mean an immersed submanifold. The second fundamental form of $N$ is defined by
\[
\mbox{\Romannum{2}}(X,Y)=\pi_{\nu N}\left(\bar{\nabla}_{\mathrm dF_0(X)}\mathrm dF_0(Y)\right),\;X,Y\in\Gamma(TN),
\]
where $\bar{\nabla}$ denotes the Levi-Civita connection of $M$ and $\pi_{\nu N}$ the orthogonal projection onto the normal bundle $\nu N$ of $N$. The mean curvature vector field $H\in\Gamma(\nu N)$ of $N$ is defined as the trace of the second fundamental form with respect to the induced Riemannian metric on $N$.

\begin{Def}
A smooth one parameter family $\{F(\cdot,t)\}_{t\in[0,T)}$, $T>0$, of immersions of $N$ into $M$ is evolving by mean curvature flow if
\begin{equation}\label{MCF}
\begin{split}
&\frac{\partial F}{\partial t}(x,t)=H(x,t),\;(x,t)\in N\times(0,T)\\&F(x,0)=F_0(x),\;x\in N.
\end{split}
\end{equation}
\end{Def}

The mean curvature flow is a quasilinear parabolic system and hence, if $N$ is compact, short time existence and uniqueness for given initial data is guaranteed by standard theory of quasilinear parabolic PDEs, see for instance Lady\v zhenskaja et al. \cite{LSU}.

From now on and throughout this paper we let $(M,J,\bar{\omega},\bar{g})$ denote a compact K\"ahler manifold of real dimension $2n$ with complex structure $J$, K\"ahler form $\bar{\omega}$, and K\"ahler metric $\bar{g}$. The K\"ahler form and K\"ahler metric are related by $\bar{g}(JX,Y)=\bar{\omega}(X,Y)$, for $X,Y\in\Gamma(TM)$. The Levi-Civita connection of $\bar{g}$ is denoted by $\bar{\nabla}$ and the Riemann curvature tensor $\bar{R}$ of $\bar{g}$ is $\bar{R}(X,Y)Z=\bar{\nabla}_X\bar{\nabla}_YZ-\bar{\nabla}_Y\bar{\nabla}_XZ-\bar{\nabla}_{[X,Y]}Z$,
for $X,Y,Z\in\Gamma(TM)$. Moreover the Ricci tensor $\bar{R}ic$ of $\bar{g}$ is $\bar{R}ic(X,Y)=\trace\;\bar{R}(.,X)Y$, for $X,Y\in\Gamma(TM)$, and the Ricci form $\bar{\rho}$, which is a real $(1,1)$-form, is defined by $\bar{\rho}(X,Y)=\bar{R}ic(JX,Y)$.

Let $L$ be a compact manifold of real dimension $n$ and $F_0:L\rightarrow M$ an immersion of $L$ into $M$. The induced Riemannian metric on $L$ is $g=F_0^*(\bar{g})$, and we set $\omega=F_0^*(\bar{\omega})$. Assume now that $F_0$ is a Lagrangian immersion, i.e. $\omega=0$. We recall some basic geometric properties of Lagrangian submanifolds. For any normal vector field $\xi\in\Gamma(\nu L)$ there is a corresponding one form $\alpha_{\xi}$ on $L$ given by $\alpha_{\xi}=F_0^*(\xi\;\lrcorner\;\bar{\omega})$. The one-form $\alpha_H=F_0^*(H\;\lrcorner\;\omega)$ is called the mean curvature form and it satisfies the following important relation first proved by Dazord in \cite{Dazord}:

\begin{Prop}\label{RicciForm}
The mean curvature form $\alpha_H$ satisfies
\[
\mathrm d\alpha_H=F_0^*(\bar{\rho}).
\]
\end{Prop}

In particular by Cartan's formula we find
\[
F_0^*\left(\mathcal{L}_H\bar{\omega}\right)=F_0^*\left(\mathrm d\left(H\;\lrcorner\;\bar{\omega}\right)\right)+F_0^*\left(H\;\lrcorner\;\mathrm d\bar{\omega}\right)=F_0^*(\bar{\rho}).
\]
Hence, if $M$ is K\"ahler-Einstein, i.e. $\bar{\rho}=\lambda\bar{\omega}$ for some $\lambda\in\mathbb{R}$, then it follows that the deformation of a Lagrangian submanifold in the direction of the mean curvature vector field is an infinitesimal symplectic motion. A natural question that arises now is whether the Lagrangian condition is preserved under the mean curvature flow. This question was answered positively by Smoczyk in \cite{Smoczyk}:

\begin{Thm}
Let $L$ be a compact $n$-dimensional manifold and let $F_0:L\rightarrow M$ be a Lagrangian immersion into a compact K\"ahler-Einstein manifold $M$. Then the mean curvature flow admits a unique smooth solution for a short time and this solution consists of Lagrangian submanifolds.
\end{Thm}

\section{Generalized Lagrangian mean curvature flow in K\"ahler manifolds that are almost Einstein}

\begin{Def}
An $n$-dimensional K\"ahler manifold $(M,J,\bar{\omega},\bar{g})$ is called almost Einstein if
\[
\bar{\rho}=\lambda\bar{\omega}+n\mathrm d\mathrm d^c\psi
\]
for some constant $\lambda\in\mathbb{R}$ and some smooth function $\psi$ on $M$.
\end{Def}

From now on we additionally assume that our K\"ahler manifold $(M,J,\bar{\omega},\bar{g})$ is almost Einstein. Given an immersion $F_0:N\rightarrow M$ of a manifold $N$ into $M$ we define a normal vector field $K\in\Gamma(\nu N)$ along $N$ by
\[
K=H-n\pi_{\nu N}\left(\bar{\nabla}\psi\right).
\]
We call $K$ the generalized mean curvature vector field of $N$. Now let $L$ be an $n$-dimensional manifold and $F_0:L\rightarrow M$ a Lagrangian immersion. Then the deformation of $L$ in direction of the generalized mean curvature vector field is an infinitesimal symplectic motion. Indeed by Dazord's result we have
\[
F_0^*\left(\mathcal{L}_K\bar{\omega}\right)=\mathrm d\alpha_H+nF_0^*\left(\mathrm d\left(\mathrm d\psi\circ J\right)\right)=F_0^*\left(\bar{\rho}-n\mathrm d\mathrm d^c\psi\right)=\lambda F_0^*(\bar{\omega})=0.
\]
Also observe that if $M$ is K\"ahler-Einstein, then $K$ is the mean curvature vector field.

In the remainder we study the generalized mean curvature flow
\begin{equation}\label{GMCF}
\begin{split}
&\frac{\partial F}{\partial t}(x,t)=K(x,t),\;(x,t)\in L\times(0,T)\\&F(x,0)=F_0(x),\;x\in L,
\end{split}
\end{equation}
for a given Lagrangian immersion $F_0:L\rightarrow M$ of a compact $n$-dimensional manifold $L$ into $M$ and $\{F(\cdot,t)\}_{t\in[0,T)}$ a smooth one-parameter family of immersions of $L$ into $M$. In order to establish the short time existence and uniqueness of this flow observe that $K$ as a differential operator differs from $H$ just by lower order terms. Hence $K$ and $H$ have the same principal symbol, so short time existence and uniqueness for (\ref{GMCF}) follows immediately.

Now let $\{F(\cdot,t)\}_{t\in[0,T)}$ be the solution to the generalized mean curvature flow (\ref{GMCF}). In the remaining part of this chapter we show that $F(\cdot,t):L\rightarrow M$ is Lagrangian for each $t\in(0,T)$. As before we denote $g=F(\cdot,t)^*(\bar{g})$ and $\omega=F(\cdot,t)^*(\bar{\omega})$. Furthermore $\nabla$ will denote the Levi-Civita connection of $g$ and $R$ the Riemannian curvature tensor of $g$. Let $p\in L$ and choose normal coordinates $\{x^i\}$ on $L$ around $p$ at time $t\in(0,T)$ and coordinates $\{y^{\alpha}\}$ on $M$ around $F(p,t)$. We have to introduce some notation. We denote $e_i=\frac{\partial F}{\partial x^i}(\cdot,t)$ and we define tensors $N$ and $\eta$ by $N_i=N(e_i)=\pi_{\nu L}\left(Je_i\right)$ and $\eta_{ij}=\eta(e_i,e_j)=\bar{g}\left(Ne_i,Ne_j\right)$. Moreover we set $h_{ijk}=h(e_i,e_j,e_k)=-\bar{g}\left(Ne_i,\bar{\nabla}_{e_j}e_k\right)$. Observe that $h_{ijk}$ is symmetric in the last two indices and fully symmetric if $F(\cdot,t)$ is Lagrangian. We also denote $\bar{R}_{klj\underline{i}}=\bar{R}(e_k,e_l,e_j,N(e_i))$. The following formula proved by Smoczyk \cite[Lem. 1.4]{Smoczyk} will be of use later:

\begin{Lem}\label{formula}
\begin{eqnarray*} &&\nabla_lh_{kij}-\nabla_kh_{lij}=\bar{R}_{klj\underline{i}}+\nabla_j\nabla_i\omega_{lk}+\omega_i^{\;m}\bar{R}_{kljm}+\omega_k^{\;m}R_{ljim}\\&&\;\;\;\;\;\;\;\;\;\;\;\;\;\;\;\;\;\;\;\;\;\;\;\;\;\;\;\;\;\;\;\;\;\;\;\;\;\;\;\;\;\;\;\;\;\;\;\;\;\;+\omega_l^{\;m}R_{jkim}+\eta^{mn}\omega_n^{\;s}\left(h_{mlj}h_{ski}-h_{mkj}h_{sli}\right).
\end{eqnarray*}
\end{Lem}

We start by computing the evolution equations of $g_{ij}$ and $\omega_{ij}$ at $p\in L$ and time $t$.

\begin{Lem}\label{EvolutionEquations}
\begin{eqnarray*}
&&i)\;\frac{\mathrm d}{\mathrm dt}\omega_{ij}=\left(\mathrm d\alpha_K\right)_{ij}\;\;\;\;\;\;\;\;\;\;\;\;\;\;\;\;\;\;\;\;\;\;\;\;\;\;\;\;\;\;\;\;\;\;\;\;\;\;\;\;\;\;\;\;\;\;\;\;\;\;\;\;\;\;\;\;\;\;\;\;\;\;\;\;\;\;\;\;\;\;\;\;\;\;\;\;\;\;\;\;\;\;\;\;\;\;\;\;\;\;\;\;\;\;\\
&&ii)\;\frac{\mathrm d}{\mathrm dt}g_{ij}=-2\eta^{mn}(\alpha_H)_mh_{nij}+2n\mathrm d\psi(\mbox{\Romannum{2}}_{ij}).
\end{eqnarray*}
\end{Lem}

\begin{proof}
\begin{eqnarray*}
\frac{\mathrm d}{\mathrm dt}\omega_{ij}&=&\frac{\mathrm d}{\mathrm dt}\bar{\omega}_{\alpha\beta}\frac{\partial F^{\alpha}}{\partial x^i}\frac{\partial F^{\beta}}{\partial x^j}\\&=&\bar{\omega}_{\alpha\beta}\left\{\frac{\partial}{\partial x^i}\frac{\partial F^{\alpha}}{\partial t}\frac{\partial F^{\beta}}{\partial x^j}+\frac{\partial F^{\alpha}}{\partial x^i}\frac{\partial}{\partial x^j}\frac{\partial F^{\beta}}{\partial t}\right\}\\&=&\bar{\omega}_{\alpha\beta}\left\{\frac{\partial K^{\alpha}}{\partial x^i}\frac{\partial F^{\beta}}{\partial x^j}+\frac{\partial F^{\alpha}}{\partial x^i}\frac{\partial K^{\beta}}{\partial x^j}\right\}\\&=&\bar{\omega}\left(\frac{\partial K}{\partial x^i},\frac{\partial F}{\partial x^j}\right)-\bar{\omega}\left(\frac{\partial K}{\partial x^j},\frac{\partial F}{\partial x^i}\right)=\left(\mathrm d\alpha_K\right)_{ij}.
\end{eqnarray*}
\begin{eqnarray*}
\frac{\mathrm d}{\mathrm dt}g_{ij}&=&\bar{g}\left(\frac{\partial}{\partial x^i}\frac{\partial F}{\partial t},\frac{\partial F}{\partial x^j}\right)+\bar{g}\left(\frac{\partial F}{\partial x^i},\frac{\partial}{\partial x^j}\frac{\partial F}{\partial t}\right)\\&=&\bar{g}\left(\frac{\partial K}{\partial x^i},\frac{\partial F}{\partial x^j}\right)+\bar{g}\left(\frac{\partial K}{\partial x^j},\frac{\partial F}{\partial x^i}\right)\\&=&-2\bar{g}\left(H,\frac{\partial^2F}{\partial x^i\partial x^j}\right)+2n\bar{g}\left(\bar{\nabla}\psi,\pi_{\nu L}\left(\frac{\partial^2F}{\partial x^i\partial x^j}\right)\right)\\&=&-2\eta^{mn}(\alpha_H)_mh_{nij}+2n\mathrm d\psi\left(\mbox{\Romannum{2}}_{ij}\right).
\end{eqnarray*}
\end{proof}

Using Lemma \ref{formula} and Lemma \ref{EvolutionEquations} we can now proceed as in \cite{Smoczyk} to prove the following lemma.

\begin{Lem}\label{EvolutionInequality}
Let $0<\tau<T$, then there exists a constant $C>0$ such that for all $t\in[0,\tau]$
\[
\frac{\mathrm d}{\mathrm dt}|\omega|^2\leq\Delta|\omega|^2+n\mathrm d\psi\left(\nabla|\omega|^2\right)+C|\omega|^2.
\]
\end{Lem}

\begin{proof}
Denote $Y=\pi_{\nu L}(\bar{\nabla}\psi)$, so that $K=H-nY$. Then
\allowdisplaybreaks
\begin{eqnarray*}
\frac{\mathrm d}{\mathrm dt}|\omega|^2&=&\frac{\mathrm d}{\mathrm dt}g^{mk}g^{jl}\omega_{mj}\omega_{kl}\\&=&-2\omega^{kl}\omega^m_{\;\;\;l}\frac{\mathrm d}{\mathrm dt}g_{mk}+2\omega^{kl}\frac{\mathrm d}{\mathrm dt}\omega_{kl}\\&=&-2\omega^{kl}\omega^m_{\;\;\;l}\left(-2\eta^{st}(\alpha_H)_sh_{tmk}+2n\mathrm d\psi\left(\mbox{\Romannum{2}}_{mk}\right)\right)\\&&\;\;\;\;\;\;\;\;\;\;\;\;\;\;\;\;\;\;\;\;+2\omega^{kl}\left(\nabla_k(\alpha_H)_l-\nabla_l(\alpha_H)_k-n(\mathrm d\alpha_Y)_{kl}\right)\\&=&4\omega^{kl}\omega^m_{\;\;\;l}\eta^{st}(\alpha_H)_sh_{tmk}-4n\omega^{kl}\omega^m_{\;\;\;l}\mathrm d\psi\left(\mbox{\Romannum{2}}_{mk}\right)\\&&\;\;\;\;\;\;\;\;\;\;\;\;\;\;\;\;\;\;\;\;+2\omega^{kl}\left(\nabla_k(\alpha_H)_l-\nabla_l(\alpha_H)_k\right)-2n\omega^{kl}(\mathrm d\alpha_Y)_{kl}\\&=&4\omega^{kl}\omega^m_{\;\;\;l}\eta^{st}(\alpha_H)_sh_{tmk}-4n\omega^{kl}\omega^m_{\;\;\;l}d\psi\left(\mbox{\Romannum{2}}_{mk}\right)\\&&\;\;\;\;\;\;\;\;\;\;\;\;\;\;\;\;\;\;\;\;+2\omega^{kl}g^{pq}\left(\nabla_kh_{lpq}-\nabla_lh_{kpq}\right)-2n\omega^{kl}(\mathrm d\alpha_Y)_{kl}\\&=&4\omega^{kl}\omega^m_{\;\;\;l}\eta^{st}(\alpha_H)_sh_{tmk}-4n\omega^{kl}\omega^m_{\;\;\;l}\mathrm d\psi\left(\mbox{\Romannum{2}}_{mk}\right)\\&&\;\;\;\;\;+2\omega^{kl}g^{pq}\left(\bar{R}_{lkq\underline{p}}+\nabla_q\nabla_p\omega_{kl}+\omega_p^{\;s}\bar{R}_{lkqs}+\omega_l^{\;s}R_{kqps}+\omega_k^{\;s}R_{qlps}\right.\\&&\;\;\;\;\;\;\;\;\;\;\;\;\;\;\;\;\;\;\;\;+\eta^{mt}\omega_t^{\;\;s}\left(h_{mkq}h_{slp}-h_{mlq}h_{skp}\right)\Big)-2n\omega^{kl}(d\alpha_Y)_{kl}\\&=&4\omega^{kl}\omega^m_{\;\;\;l}\eta^{st}(\alpha_H)_sh_{tmk}-4n\omega^{kl}\omega^m_{\;\;\;l}\mathrm d\psi\left(\mbox{\Romannum{2}}_{mk}\right)+2\omega^{kl}\bar{R}_{lk\;\underline{p}}^{\;\;\;p}+\Delta|\omega|^2\\&&\;\;\;\;\;-|\nabla\omega|^2+2\omega^{kl}\omega_p^{\;\;s}\bar{R}_{lk\;s}^{\;\;\;p}+2\omega^{kl}\omega_l^{\;\;s}R_{k\;ps}^{\;\;p}+2\omega^{kl}\omega_k^{\;\;s}R_{\;lps}^{p}\\&&\;\;\;\;\;\;\;\;\;\;\;\;\;\;\;\;\;\;\;\;+2\omega^{kl}\eta^{mt}\omega_t^{\;\;s}\left(h_{mk}^{\;\;\;\;\;p}h_{slp}-h_{ml}^{\;\;\;\;\;p}h_{skp}\right)-2n\omega^{kl}(\mathrm d\alpha_Y)_{kl}.
\end{eqnarray*}
For terms of the form $\omega^{sl}\omega^m_{\;\;\;i}T_{\;slm}^i$ we have
\begin{eqnarray*}
2\omega^{sl}\omega^m_{\;\;\;i}T_{\;slm}^i&=&2\sum_{s,l}\left(\omega_{sl}\sum_{m,i}\omega_{mi}T_{islm}\right)\\&\leq&\sum_{s,l}(\omega_{sl})^2+\sum_{s,l}\left(\sum_{m,i}\omega_{ml}T_{islm}\right)^2\\&\leq&|\omega|^2+n^2\sum_{s,l,m,i}(\omega_{ml})^2(T_{islm})^2\leq\left(1+n^2|T|^2\right)|\omega|^2.
\end{eqnarray*}
Since $L\times[0,\tau]$ is compact we can choose a constant $C>0$ such that for all $t\in[0,\tau]$
\begin{eqnarray*}
\frac{\mathrm d}{\mathrm dt}|\omega|^2&\leq&\Delta|\omega|^2+C|\omega|^2+2\omega^{kl}\bar{R}_{lk\;\underline{p}}^{\;\;\;p}-2n\omega^{kl}(\mathrm d\alpha_Y)_{lk}.
\end{eqnarray*}
It remains to find an estimate for the last two terms. We have
\begin{eqnarray*}
(\alpha_Y)_l&=&\bar{\omega}\left(\pi_{\nu L}\left(\bar{\nabla}\psi\right),e_l\right)\\&=&-\bar{g}\left(\bar{\nabla}\psi,J(e_l)\right)-g^{mk}\bar{g}\left(\bar{\nabla}\psi,e_k\right)\bar{\omega}\left(e_m,e_l\right)\\&=&\mathrm d^c\psi(e_l)-g^{mk}\mathrm d\psi(e_k)\omega_{ml},
\end{eqnarray*}
and hence
\begin{eqnarray*}
(\mathrm d\alpha_Y)_{kl}&=&\frac{\partial}{\partial x^k}\mathrm d^c\psi(e_l)-\frac{\partial}{\partial x^l}\mathrm d^c\psi(e_k)\\&&\;\;\;\;\;\;\;\;\;\;\;\;-\frac{\partial}{\partial x^k}\left(g^{mj}\mathrm d\psi(e_j)\omega_{ml}\right)+\frac{\partial}{\partial x^l}\left(g^{mj}\mathrm d\psi(e_j)\omega_{mk}\right)\\&=&\mathrm d\mathrm d^c\psi\left(e_k,e_l\right)-g^{mj}\omega_{ml}\frac{\partial}{\partial x^k}\mathrm d\psi(e_j)\\&&\;\;\;\;\;+g^{mj}\omega_{mk}\frac{\partial}{\partial x^l}\mathrm d\psi(e_j)-g^{mj}\mathrm d\psi(e_j)\frac{\partial\omega_{ml}}{\partial x^k}+g^{mj}\mathrm d\psi(e_j)\frac{\partial\omega_{mk}}{\partial x^l}.
\end{eqnarray*}
Multiplying both sides with $-2n\omega^{kl}$, using the K\"ahler and the almost Einstein condition, and estimating the quadratic terms in $\omega$ we get
\begin{eqnarray*}
-2n\omega^{kl}(\mathrm d\alpha_Y)_{kl}\leq C|\omega|^2-2\omega^{kl}\bar{\rho}_{kl}+2ng^{mj}\omega^{kl}\mathrm d\psi(e_j)\left(\frac{\partial\omega_{ml}}{\partial x^k}-\frac{\partial\omega_{mk}}{\partial x^l}\right).
\end{eqnarray*}
Using that $\omega$ is closed we find
\begin{eqnarray*}
&&g^{mj}\omega^{kl}\left(\frac{\partial\omega_{ml}}{\partial x^k}-\frac{\partial\omega_{mk}}{\partial x^l}\right)\mathrm d\psi(e_j)=\mathrm d\psi\left(g^{mj}g^{ik}g^{sl}\omega_{is}\frac{\partial\omega_{kl}}{\partial x^m}e_j\right)\\&&\;\;\;\;\;\;\;\;\;\;\;\;\;\;\;=\mathrm d\psi\left(\frac{1}{2}g^{mj}\frac{\partial|\omega|^2}{\partial x^m}e_j\right)=\frac{1}{2}\mathrm d\psi\left(\nabla|\omega|^2\right).
\end{eqnarray*}
Putting all together yields
\[
\frac{\mathrm d}{\mathrm dt}|\omega|^2\leq\Delta|\omega|^2+n\mathrm d\psi\left(\nabla|\omega|^2\right)+C|\omega|^2+2\omega^{kl}\bar{R}_{lk\;\underline{p}}^{\;\;\;p}-2\omega^{kl}\bar{\rho}_{kl}.
\]
Now by definition of the tensor $N$ we have $N(e_p)=J(e_p)-\omega_p^{\;\;m}e_m$ and so
\[
\bar{R}_{lkq\underline{p}}=\bar{R}(e_l,e_k,e_q,J(e_p))-\omega_p^{\;\;m}\bar{R}_{lkqm}.
\]
Multiplying both sides with $2\omega^{kl}g^{pq}$ and estimating the quadratic term in $\omega$ gives
\[
2\omega^{kl}\bar{R}_{lk\;\underline{p}}^{\;\;\;p}\leq2\omega^{kl}g^{pq}\bar{R}(e_l,e_k,e_q,J(e_p))+C|\omega|^2.
\]
Using the following well known identity from K\"ahler geometry
\[
g^{pq}\bar{R}(e_l,e_k,e_q,J(e_p))=\bar{\rho}_{kl},
\]
we finally obtain
\[
\frac{\mathrm d}{\mathrm dt}|\omega|^2\leq\Delta|\omega|^2+n\mathrm d\psi\left(\nabla|\omega|^2\right)+C|\omega|^2.
\]
\end{proof}

Applying the parabolic maximum principle we conclude that $F(\cdot,t):L\rightarrow M$ is Lagrangian for each $t\in[0,T)$. This motivates the following definition:
\begin{Def}
A family of Lagrangian submanifolds satisfying (\ref{GMCF}) is said to evolve by generalized Lagrangian mean curvature flow.
\end{Def}

And we have proved the following theorem:

\begin{Thm}\label{MainTheorem}
Let $L$ be a compact $n$-dimensional manifold and $F_0:L\rightarrow M$ a Lagrangian immersion of $L$ into a compact K\"ahler manifold $M$ that is almost Einstein. Then the generalized Lagrangian mean curvature flow admits a unique smooth solution for a short time and this solution consists of Lagrangian submanifolds.
\end{Thm}

\section{A variational approach to the generalized mean curvature flow}
Let $\mathcal{S}$ be the infinite dimensional manifold consisting of all compact $n$-dimensional submanifolds of $M$. In this chapter we show that the generalized mean curvature flow is the gradient flow of a volume functional on $\mathcal{S}$. Let $N\in\mathcal{S}$, then the tangent space of $\mathcal{S}$ at $N$ consists of the normal vector fields along $N$ and for any Riemannian metric $g$ on $M$ there is a natural $L^2$-metric on $\mathcal{S}$ given by
\[
\langle Y,Z\rangle_{g,L^2}=\int_Ng(Y,Z)\mathrm dV_g,
\]
for $Y,Z\in\Gamma(\nu N)$.

We define two conformally rescaled Riemannian metrics $\tilde{g}$ and $\hat{g}$ on $M$ by
\[
\tilde{g}=e^{2\psi}\bar{g}\quad\mbox{and}\quad\hat{g}=e^{\frac{2n}{n+2}\psi}\bar{g}.
\]
Then we have the following variational characterization of the generalized mean curvature flow:

\begin{Prop}\label{GradientFlow}
The generalized mean curvature flow is the negative gradient flow of the volume functional $Vol_{\tilde{g}}$ on $\mathcal{S}$ with respect to the $L^2$-metric $\langle\cdot,\cdot\rangle_{\hat{g},L^2}$.
\end{Prop}

\begin{proof}
Let $N\in\mathcal{S}$ and let $Y$ be a normal vector field along $N$. Then the first variation of the volume functional gives
\[
\delta_Y\Vol_{\tilde{g}}(N)=-\int_N\tilde{g}(Y,\tilde{H})\mathrm dV_{\tilde{g}},
\]
where $\tilde{H}$ is the mean curvature vector field on $N$ with respect to the metric on $N$ which is induces by $\tilde{g}$. It is easy to show that
\[
\tilde{H}=e^{-2\psi}\left(H-n\pi_{\nu N}\left(\bar{\nabla}\psi\right)\right).
\]
Hence
\begin{eqnarray*}
\delta_Y\Vol_{\tilde{g}}(N)&=&-\int_Ne^{n\psi}\bar{g}(H-n\pi_{\nu N}(\bar{\nabla}\psi),Y)\mathrm dV_{\bar{g}}\\&=&-\int_Ne^{\left(n-\frac{2n}{n+2}-\frac{2n}{n+2}\frac{n}{2}\right)\psi}\hat{g}(K,Y)\mathrm dV_{\hat{g}}\\&=&-\int_N\hat{g}(K,Y)\mathrm dV_{\hat{g}}=-\langle K,Y\rangle_{\hat{g},L^2}.
\end{eqnarray*}
\end{proof}

\section{The case of almost Calabi-Yau manifolds}
We introduce almost Calabi-Yau manifolds and special Lagrangian submanifolds as defined by Joyce in \cite[\S8.4]{JoyceBook}.

\begin{Def}
An $n$-dimensional almost Calabi-Yau manifold $(M,J,\bar{\omega},\bar{g},\Omega)$ is an $n$-dimensional K\"ahler manifold $(M,J,\bar{\omega},\bar{g})$ together with a non-vanishing holomorphic volume form $\Omega$.
\end{Def}

Given an $n$-dimensional almost Calabi-Yau manifold $(M,J,\bar{\omega},\bar{g},\Omega)$ we can define a smooth function $\psi$ on $M$ by
\[
e^{2n\psi}\frac{\bar{\omega}^n}{n!}=\left(-1\right)^{\frac{n(n-1)}{2}}\left(\frac{i}{2}\right)^n\Omega\wedge\bar{\Omega}.
\]
Here $\bar{\Omega}$ denotes the complex conjugate of $\Omega$. Then $(M,J,\bar{\omega},\bar{g},\Omega)$ is Calabi-Yau if and only if $\psi$ vanishes identically. Using $|\Omega|_{\bar{g}}=2^{\frac{n}{2}}e^{n\psi}$ and the following formula for the Ricci form of a K\"ahler manifold with trivial canonical bundle (see for instance \cite[\S 7.1]{JoyceBook})
\[
\bar{\rho}=\mathrm d\mathrm d^c\log|\Omega|_{\bar{g}}
\]
we find
\[
\bar{\rho}=n\mathrm d\mathrm d^c\psi.
\]
Hence almost Calabi-Yau manifolds are almost Einstein and Theorem \ref{MainTheorem} holds in this case. Let $\tilde{g}$ be a conformally rescaled metric on $M$ defined by $\tilde{g}=e^{2\psi}\bar{g}$. One easily proves that $\realpart\;\Omega$ is a calibrating $n$-form on $(M,\tilde{g})$. This leads to the definition of special Lagrangian submanifolds in almost Calabi-Yau manifolds.

\begin{Def}
An oriented Lagrangian submanifold $L$ of an almost Calabi-Yau manifold $M$ is called special Lagrangian if it is calibrated with respect to $\realpart\;\Omega$ for the metric $\tilde{g}$. More generally, an oriented Lagrangian submanifold $L$ is special Lagrangian with phase $e^{i\theta_0}\in\mathbb{R}$, if $L$ is calibrated with respect to $\realpart(e^{-i\theta_0}\Omega)$ for the metric $\tilde{g}$.
\end{Def}

Besides the fact that one is able to write down explicit examples of almost Calabi-Yau metrics on compact manifolds there is another reason for studying almost Calabi-Yau manifolds. Recall that by the theorem of Tian and Todorov the moduli space $\mathcal{M}_{CY}$ of Calabi-Yau metrics of a compact Calabi-Yau manifold is of dimension $h^{1,1}(M)+2h^{n-1,1}(M)+1$, where $h^{i,j}(M)$ are the Hodge numbers of $M$. In particular $\mathcal{M}_{CY}$ is finite dimensional. In the study of moduli spaces of $J$-holomorphic curves in symplectic manifolds it turns out that for a generic almost complex structure $J$ the moduli space $\mathcal{M}_J$ of embedded $J$-holomorphic curves is a smooth manifold, while for a fixed almost complex structure $J$ the space $\mathcal{M}_J$ can have singularities (see \cite{McDuffSalamon} for details). Now the moduli space $\mathcal{M}_{ACY}$ of almost Calabi-Yau structures is of infinite dimension and therefore choosing a generic almost Calabi-Yau metric is a more powerful thing to do than choosing a generic Calabi-Yau metric. We explain why this is of certain interest. It was proved by McLean \cite{McLean} that the moduli space of compact special Lagrangian submanifolds $\mathcal{M}_{SL}$ in a Calabi-Yau manifold is a smooth manifold of dimension $b^1(L)$, the first Betti number of $L$. An important question is whether it is possible to compactify $\mathcal{M}_{SL}$ in order to define invariants of Calabi-Yau manifolds by counting special Lagrangian submanifolds. One approach to this problem, due to Joyce, is to study the moduli space of special Lagrangian submanifolds with conical singularities in almost Calabi-Yau manifolds (see \cite{Joyce5} for a survey of his results). In particular Joyce conjectures that for generic almost Calabi-Yau metrics the moduli space of special Lagrangian submanifolds with conical singularities is a smooth finite dimensional manifold.

We come back to the study of the generalized Lagrangian mean curvature flow. First observe that special Lagrangian submanifolds in an almost Calabi-Yau manifold $M$ are minimal with respect to $\tilde{g}$. By Proposition \ref{GradientFlow} the generalized Lagrangian mean curvature flow decreases volume with respect to $\tilde{g}$. Therefore the generalized Lagrangian mean curvature flow is in this sense the right flow to consider. Harvey and Lawson show in \cite{HarveyLawson} that
\[
F_0^*(\Omega)=e^{i\theta+n\psi}\mathrm dV_g,
\]
for $F_0:L\rightarrow M$ a Lagrangian immersion. The map $\theta:L\rightarrow S^1$ is called the Lagrangian angle of $L$. From this we obtain an alternative characterization of special Lagrangian submanifolds.

\begin{Prop}
An oriented Lagrangian submanifold $L$ is special Lagrangian with phase $e^{i\theta_0}$ if and only if
\[
(\cos{\theta_0}\;\imagpart\;\Omega-\sin{\theta_0}\;\realpart\;\Omega)|_L=0.
\]
In particular, an oriented Lagrangian submanifold is special Lagrangian with phase $e^{i\theta_0}$ if and only if the Lagrangian angle is constant.
\end{Prop}

The Lagrangian angle is closely related to the generalized Lagrangian mean curvature flow as proved in the next proposition.

\begin{Prop}\label{theta}
Let $L$ be a Lagrangian submanifold of $M$. Then
\[
\alpha_K=-\mathrm d\theta.
\]
\end{Prop}

\begin{proof}
The decomposition
\[
\Lambda^nT^*M\otimes\mathbb{C}=\bigoplus_{p+q=n}\Lambda^{p,q}T^*M
\]
is invariant under the holonomy representation of $\bar{g}$. Hence there exists a complex one form $\eta$ on $M$ satisfying $\bar{\nabla}\Omega=\eta\otimes\Omega$. Moreover, since $\Omega$ is holomorphic, $\eta$ is in fact a one form of type $(1,0)$. Using $\Omega\wedge\bar{\Omega}=e^{2n\psi}\mathrm dV_{\bar{g}}$ we find by computing $\bar{\nabla}\left(\Omega\wedge\bar{\Omega}\right)$ the equality
\[
(\eta+\bar{\eta})\otimes\Omega\wedge\bar{\Omega}=2n\mathrm d\psi\otimes\Omega\wedge\bar{\Omega}.
\]
Hence $\eta=2n\partial\psi$ and so $\bar{\nabla}\Omega=2n\partial\psi\otimes\Omega$. Following the computation by Thomas and Yau \cite[Lem. 2.1]{ThomasYau} we obtain
\[
\bar{\nabla}\Omega=(i\mathrm d\theta+n\mathrm d\psi+i\alpha_H)\otimes\Omega
\]
and establish the equality
\[
\alpha_H-n\mathrm d^c\psi=-\mathrm d\theta.
\]
But $\alpha_H-n\mathrm d^c\psi=\alpha_K$ and hence $\alpha_K=-\mathrm d\theta$.
\end{proof}

Now let $\{F(\cdot,t)\}_{t\in[0,T)}$ be the solution to the generalized mean curvature flow with initial condition $F_0:L\rightarrow M$ a Lagrangian immersion. Then we have the following proposition:

\begin{Prop}
Under the generalized Lagrangian mean curvature flow the Lagrangian angle of $L$ satisfies
\[
\frac{\mathrm d}{\mathrm dt}\theta=\Delta\theta+n\mathrm d\psi(\nabla\theta).
\]
\end{Prop}

\begin{proof}
On the one hand
\[
\frac{\mathrm d}{\mathrm dt}e^{i\theta+n\psi}\mathrm dV_g=i\frac{\mathrm d\theta}{\mathrm dt}e^{i\theta+n\psi}\mathrm dV_g+n\mathrm d\psi(K)e^{i\theta+n\psi}\mathrm dV_g+e^{i\theta+n\psi}\frac{\mathrm d}{\mathrm dt}\mathrm dV_g
\]
and on the other hand, using $F(\cdot,t)^*(\Omega)=e^{i\theta+n\psi}\mathrm dV_g$, we have
\[
\frac{\mathrm d}{\mathrm dt}e^{i\theta+n\psi}\mathrm dV_g=F(\cdot,t)^*(\mathcal{L}_K\Omega)=F(\cdot,t)^*(\mathrm d(K\;\lrcorner\;\Omega))+F(\cdot,t)^*(K\;\lrcorner\;\mathrm d\Omega).
\]
Since $\Omega$ is holomorphic, $\mathrm d\Omega=0$. Moreover by Proposition \ref{theta} we have $K=J(\nabla\theta)$ and hence
\begin{eqnarray*}
F(\cdot,t)^*(\mathrm d(K\;\lrcorner\;\Omega))&=&iF(\cdot,t)^*(\mathrm d(\nabla\theta\;\lrcorner\;\Omega))=i\mathrm d(e^{i\theta+n\psi}\nabla\theta\;\lrcorner\;\mathrm dV_g)\\&=&ie^{i\theta+n\psi}\left(\mathrm d(\nabla\theta\;\lrcorner\;\mathrm dV_g)+n\mathrm d\psi\wedge(\nabla\theta\;\lrcorner\;\mathrm dV_g)\right)\\&&\;\;\;\;\;\;\;\;\;\;\;-e^{i\theta+n\psi}\mathrm d\theta\wedge(\nabla\theta\;\lrcorner\;\mathrm dV_g)\\&=&ie^{i\theta+n\psi}\left(\Delta\theta+n\mathrm d\psi(\nabla\theta)\right)\mathrm dV_g-e^{i\theta+n\psi}|\nabla\theta|^2\mathrm dV_g.
\end{eqnarray*}
Comparing imaginary parts yields
\[
\frac{\mathrm d}{\mathrm dt}\theta=\Delta\theta+n\mathrm d\psi(\nabla\theta).
\]
\end{proof}

We end this paper by showing how the generalized Lagrangian mean curvature flow in an almost Calabi-Yau manifold can be integrated to a scalar equation. Let $\hat{\omega}$ be the canonical symplectic structure on the cotangent bundle $T^*L$ of $L$. Then by the Lagrangian neighbourhood theorem \cite[Thm. 3.33]{McDuffSalamon2} there exists an immersion $\Phi:U\rightarrow V$ from an open neighbourhood $U$ of the zero section in $T^*L$ onto an open neighbourhood $V$ of $L$ in $M$, such that $\hat{\omega}=\Phi^*(\bar{\omega})$ and $\Phi(x,0)=F_0(x)$ for $x\in L$. It is not hard to see that all Lagrangian submanifolds in $M$ which are $C^1$-close to $L$ correspond to graphs in $T^*L$ of closed one-forms on $L$.

\begin{Thm}
Let $F_0:L\rightarrow M$ be a zero Maslov class Lagrangian, i.e. $\theta:L\rightarrow\mathbb{R}$ is a well defined smooth function on $L$, let $\Phi$ be as above, and let $\{u(\cdot,t)\}_{t\in[0,T)}$ be a smooth one-parameter family of smooth functions on $L$ satisfying
\begin{equation}\label{Integrated}
\begin{split}
&\frac{\partial u}{\partial t}(x,t)=\theta(x,t),\;(x,t)\in L\times(0,T)\\
&u(x,0)=0,\;x\in L.
\end{split}
\end{equation}
Here $\theta(\cdot,t)$ denotes the Lagrangian angle of the Lagrangian immersion $\Phi\circ \mathrm du(\cdot,t)$ of $L$ into $M$. Choosing $T>0$ sufficiently small we can assume that the graph of $\mathrm du(\cdot,t)$ lies in $U$ for $t\in[0,T)$. Then there exists a family of diffeomorphisms $\{\varphi(.,t)\}_{t\in[0,T)}$ of $L$, such that the immersions $\{F(\cdot,t)\}_{t\in[0,T)}$ of $L$ into $M$ defined by
\[
F(x,t)=\Phi(\varphi(x,t),\mathrm du(\varphi(x,t),t)),\;x\in L,
\]
evolve by generalized Lagrangian mean curvature flow.
\end{Thm}

The proof of this theorem can be found in \cite{Smoczyk1} in the case when the ambient space is $\mathbb{C}^n$. When the ambient space is a general almost Calabi-Yau manifold the proof is analogous.

{\footnotesize\noindent LINCOLN COLLEGE, TURL STREET, OX1 3DR, OXFORD, UNITED KINGDOM}\\
{\small\verb|behrndt@maths.ox.ac.uk|}
\end{document}